\documentclass[12pt]{article}
\usepackage{amsfonts}
\usepackage{latexsym,amsmath,amssymb}
\usepackage{graphics,epstopdf}
\usepackage[pdftex]{graphicx}

\numberwithin{equation}{section}
\begin{document}

\bigskip

\bigskip

\begin{center}
{\Large \textbf{\ On $(p,q)$-analogue of Bernstein Operators (Revised)}}

\bigskip

\textbf{M. Mursaleen}, \textbf{Khursheed J. Ansari} and \textbf{Asif Khan}

Department of\ Mathematics, Aligarh Muslim University, Aligarh--202002, India%
\\[0pt]

mursaleenm@gmail.com; ansari.jkhursheed@gmail.com; asifjnu07@gmail.com \\[0pt%
]

\bigskip

\bigskip

\textbf{Abstract}
\end{center}

\parindent=8mm {\footnotesize {In the present article, we have given a corrigendum to our paper ``On $(p,q)$-analogue of Bernstein operators" published in Applied Mathematics and Computation 266 (2015) 874-882. }}\newline



\section{Construction of Revised Operators}

Mursaleen et. al \cite{mka} introduced $(p,q)$-analogue of Bernstein operators as
\begin{equation*}
B_{n,p,q}(f;x)=\sum\limits_{k=0}^{n}\left[
\begin{array}{c}
n \\
k%
\end{array}%
\right] _{p,q}x^{k}\prod\limits_{s=0}^{n-k-1}(p^{s}-q^{s}x)~~f\left( \frac{%
[k]_{p,q}}{[n]_{p,q}}\right) ,~~x\in \lbrack 0,1].\eqno(1)
\end{equation*}
But $B_{n,p,q}(1;x)\ne1$ for all $x\in[0,1]$. Hence, we are re-introducing our operators as follows:
\begin{equation*}
B_{n,p,q}(f;x)=\frac1{p^{\frac{n(n-1)}2}}\sum\limits_{k=0}^{n}\left[
\begin{array}{c}
n \\
k%
\end{array}%
\right] _{p,q}p^{\frac{k(k-1)}2}x^{k}\prod\limits_{s=0}^{n-k-1}(p^{s}-q^{s}x)~~f\left( \frac{%
[k]_{p,q}}{p^{k-n}[n]_{p,q}}\right) ,~~x\in \lbrack 0,1].\eqno(2)
\end{equation*}

Note that for $p=1$, $(p,q)$-Bernstein operators given by (2) turn out to
be $q$- Bernstein operators.

We have the following basic result:\newline

\parindent=0mm\textbf{Lemma 1.} For $x\in \lbrack 0,1],~0<q<p\leq 1$, we
have\newline
\newline
(i)~~$B_{n,p,q}(1;x)=~1$;\newline
(ii)~$B_{n,p,q}(t;x)=~x$;\newline
(iii) $B_{n,p,q}(t^{2};x)=~\frac{p^{n-1}}{[n]_{p,q}}x+\frac{%
q[n-1]_{p,q}}{[n]_{p,q}}x^{2}$.\newline

\parindent=0mm\textbf{Proof}. (i)
\begin{equation*}
B_{n,p,q}(1;x)=\frac1{p^{\frac{n(n-1)}2}}\sum\limits_{k=0}^{n}\left[
\begin{array}{c}
n \\
k%
\end{array}%
\right] _{p,q}p^{\frac{k(k-1)}2}x^{k}\prod\limits_{s=0}^{n-k-1}(p^{s}-q^{s}x)=1.
\end{equation*}

(ii)
\begin{eqnarray*}
B_{n,p,q}(t;x) &=&\frac1{p^{\frac{n(n-1)}2}}\sum\limits_{k=0}^{n}\left[
\begin{array}{c}
n \\
k%
\end{array}%
\right] _{p,q}p^{\frac{k(k-1)}2}x^{k}\prod\limits_{s=0}^{n-k-1}(p^{s}-q^{s}x)~~\frac{%
[k]_{p,q}}{p^{k-n}[n]_{p,q}}\\
&=&\frac1{p^{\frac{n(n-3)}2}}\sum\limits_{k=0}^{n-1}\left[
\begin{array}{c}
n-1 \\
k%
\end{array}%
\right] _{p,q}p^{\frac{(k+1)(k-2)}2}x^{k+1}\prod\limits_{s=0}^{n-k-2}(p^{s}-q^{s}x)\\
&=&\frac x{p^{\frac{(n-1)(n-2)}2}}\sum\limits_{k=0}^{n-1}\left[
\begin{array}{c}
n-1 \\
k%
\end{array}%
\right] _{p,q}p^{\frac{k(k-1)}2}x^{k}\prod\limits_{s=0}^{n-k-2}(p^{s}-q^{s}x)=x.
\end{eqnarray*}

(iii)
\begin{eqnarray*}
B_{n,p,q}(t^2;x) &=&\frac1{p^{\frac{n(n-1)}2}}\sum\limits_{k=0}^{n}\left[
\begin{array}{c}
n \\
k%
\end{array}%
\right] _{p,q}p^{\frac{k(k-1)}2}x^{k}\prod\limits_{s=0}^{n-k-1}(p^{s}-q^{s}x)~~\frac{%
[k]^2_{p,q}}{p^{2k-2n}[n]^2_{p,q}}\\
&=&\frac1{p^{\frac{n(n-5)}2}}\sum\limits_{k=0}^{n-1}\left[
\begin{array}{c}
n-1 \\
k%
\end{array}%
\right] _{p,q}p^{\frac{(k+1)(k-4)}2}x^{k+1}\prod\limits_{s=0}^{n-k-2}(p^{s}-q^{s}x)~\frac{[k+1]_{p,q}}{[n]_{p,q}}\\
&=&\frac1{p^{\frac{n(n-5)}2}[n]_{p,q}}\sum\limits_{k=0}^{n-1}\left[
\begin{array}{c}
n-1 \\
k%
\end{array}%
\right] _{p,q}p^{\frac{(k+1)(k-4)}2}x^{k+1}\prod\limits_{s=0}^{n-k-2}(p^{s}-q^{s}x)~(p^k+q[k]_{p,q})\\
&=&\frac1{p^{\frac{n(n-5)}2}[n]_{p,q}}\sum\limits_{k=0}^{n-1}\left[
\begin{array}{c}
n-1 \\
k%
\end{array}%
\right] _{p,q}p^{\frac{k^2-k-4}2}x^{k+1}\prod\limits_{s=0}^{n-k-2}(p^{s}-q^{s}x)\\
&&+\frac {q[n-1]_{p,q}}{p^{\frac{n(n-5)}2}[n]_{p,q}}\sum\limits_{k=0}^{n-2}\left[
\begin{array}{c}
n-2 \\
k%
\end{array}%
\right] _{p,q}p^{\frac{(k+2)(k-3)}2}x^{k+2}\prod\limits_{s=0}^{n-k-3}(p^{s}-q^{s}x)\\
&=&\frac{p^{n-1}x}{[n]_{p,q}}\frac1{p^{\frac{(n-1)(n-2)}2}}\sum\limits_{k=0}^{n-1}\left[
\begin{array}{c}
n-1 \\
k%
\end{array}%
\right] _{p,q}p^{\frac{k(k-1)}2}x^{k}\prod\limits_{s=0}^{n-k-2}(p^{s}-q^{s}x)\\
&&+\frac {q[n-1]_{p,q}~x^2}{[n]_{p,q}}\frac1{p^{\frac{(n-2)(n-3)}2}}\sum\limits_{k=0}^{n-2}\left[
\begin{array}{c}
n-2 \\
k%
\end{array}%
\right] _{p,q}p^{\frac{k(k-1)}2}x^{k}\prod\limits_{s=0}^{n-k-3}(p^{s}-q^{s}x)\\
&=&\frac{p^{n-1}}{[n]_{p,q}}x+\frac{%
q[n-1]_{p,q}}{[n]_{p,q}}x^{2}.
\end{eqnarray*}
\newline

Now, we prove Korovkin's type approximation theorem.\\

\parindent=0mm\textbf{Theorem 1}. Let $0<q_{n}<p_{n}\leq 1$ such that $%
\lim\limits_{n\rightarrow \infty }p_{n}=1$ and $\lim\limits_{n\rightarrow
\infty }q_{n}=1$. Then for each $f\in C[0,1],~B_{n,p_{n},q_{n}}(f;x)$
converges uniformly to $f$ on $[0,1]$.\\

\parindent=0mm\textbf{Proof}. By the Korovkin's Theorem it suffices to show
that
\begin{equation*}
\lim\limits_{n\rightarrow \infty }\Vert
B_{n,p_{n},q_{n}}(t^{m};x)-x^{m}\Vert _{C[0,1]}=0,~~~m=0,1,2.
\end{equation*}%
By Lemma 1(i)-(ii), it is clear that
\begin{equation*}
\lim\limits_{n\rightarrow \infty }\Vert B_{n,p_{n},q_{n}}(1;x)-1\Vert
_{C[0,1]}=0;
\end{equation*}%
\begin{equation*}
\lim\limits_{n\rightarrow \infty }\Vert B_{n,p_{n},q_{n}}(t;x)-x\Vert
_{C[0,1]}=0.
\end{equation*}%
Using $q_{n}[n-1]_{p_{n},q_{n}}=[n]_{p_{n},q_{n}}-p_{n}^{n-1}$ and by Lemma 1 (iii), we have
\begin{eqnarray*}
|B_{n,p_{n},q_{n}}(t^{2};x)-x^{2}|_{C[0,1]} &=&\biggl{|}\frac{%
p_{n}^{n-1}x}{[n]_{p_{n},q_{n}}}+\biggl{(}\frac{%
q_{n}[n-1]_{p_{n},q_{n}}}{[n]_{p_{n},q_{n}}}-1\biggl{)}x^{2}\biggl{|} \\
&\leq &\frac{%
p_{n}^{n-1}}{[n]_{p_{n},q_{n}}}x+\frac{p_{n}^{n-1}}{%
[n]_{p_{n},q_{n}}}x^{2}.
\end{eqnarray*}%
Taking maximum of both sides of the above inequality, we get\newline
\begin{equation*}
\Vert B_{n,p_{n},q_{n}}(t^{2};x)-x^{2}\Vert _{C[0,1]}\leq \frac{2p_{n}^{n-1}%
}{[n]_{p_{n},q_{n}}}
\end{equation*}%
which yields
\begin{equation*}
\lim\limits_{n\rightarrow \infty }\Vert
B_{n,p_{n},q_{n}}(t^{2};x)-x^{2}\Vert _{C[0,1]}=0.
\end{equation*}%
Thus the proof is completed.\newline

\newpage

\section{Example}

With the help of Matlab, we show comparisons and some illustrative
graphics for the convergence of operators (2) to the function $ f(x)=(x-\frac{1}{3})(x-\frac{1}{2})(x-\frac{3}{4})$ under different parameters.

From figure (1) we can observe that as the value the $q$ increases, $(p,q)$-Bernstein operators given by (2) converges towards the function.
\begin{figure*}[htb!]
\begin{center}
\includegraphics[height=6cm, width=10cm]{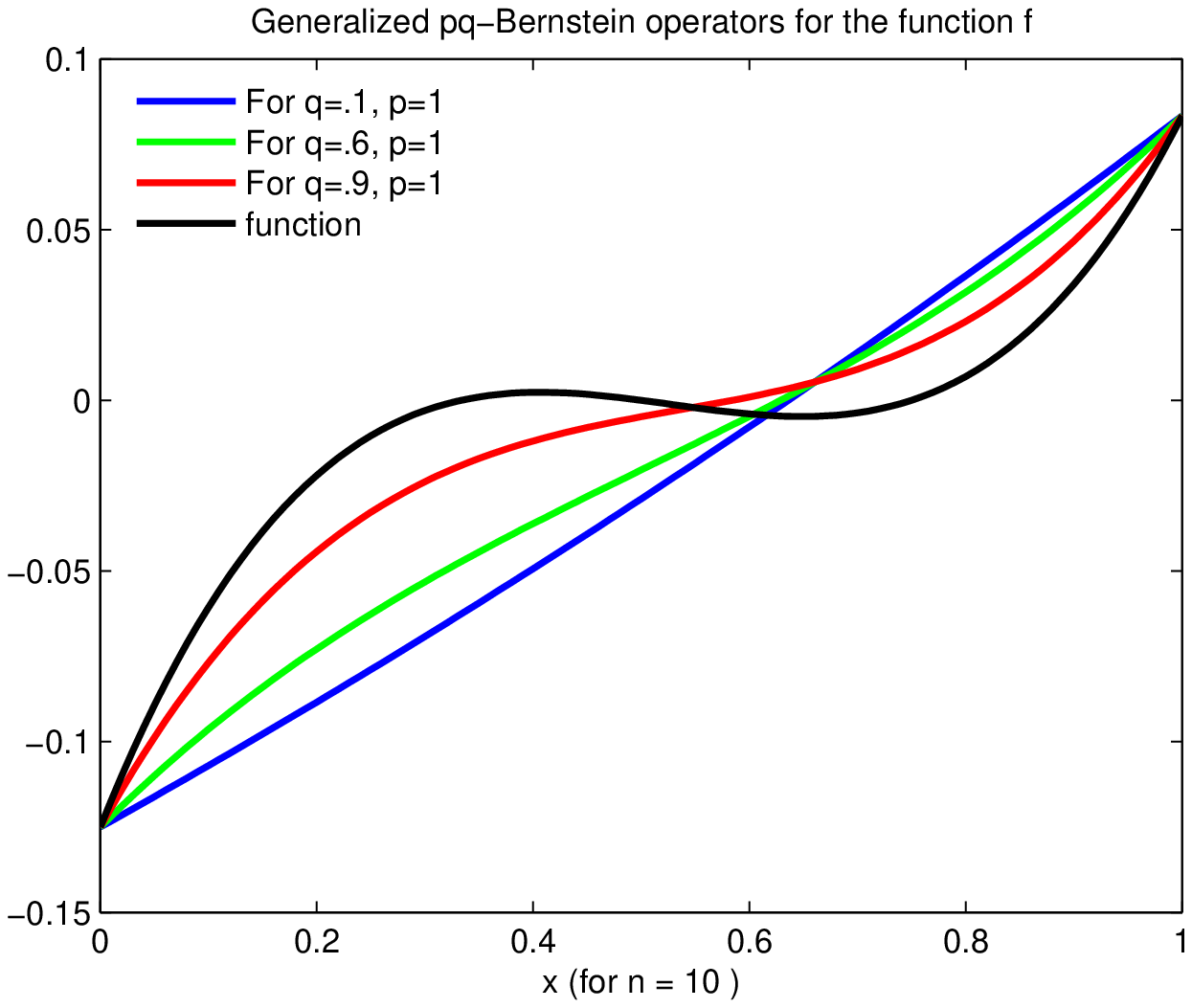}
\end{center}
\caption{}
\end{figure*}

\begin{figure*}[htb!]
\begin{center}
\includegraphics[height=6cm, width=10cm]{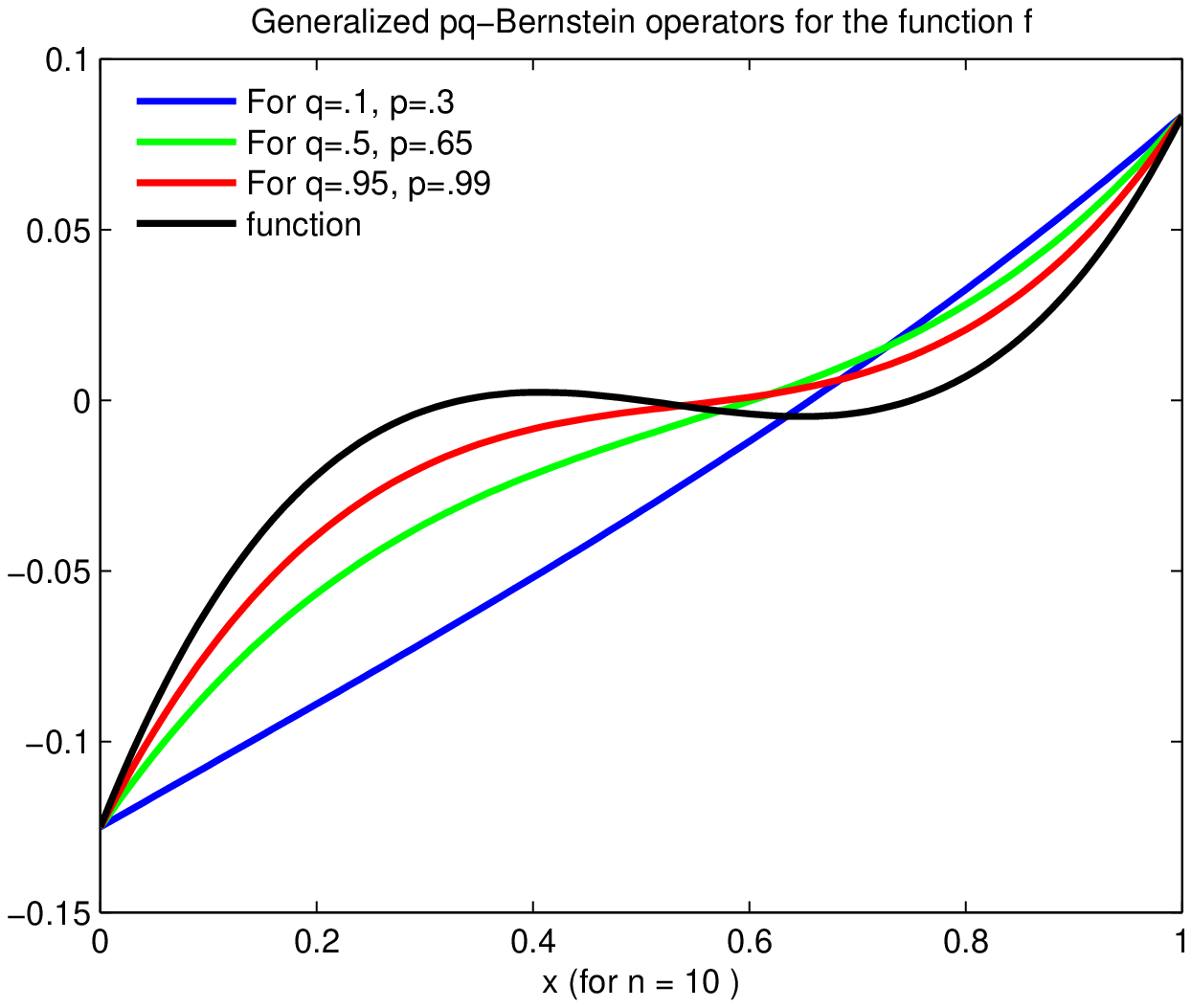}
\end{center}
\caption{}
\end{figure*}

\begin{figure*}[htb!]
\begin{center}
\includegraphics[height=6cm, width=10cm]{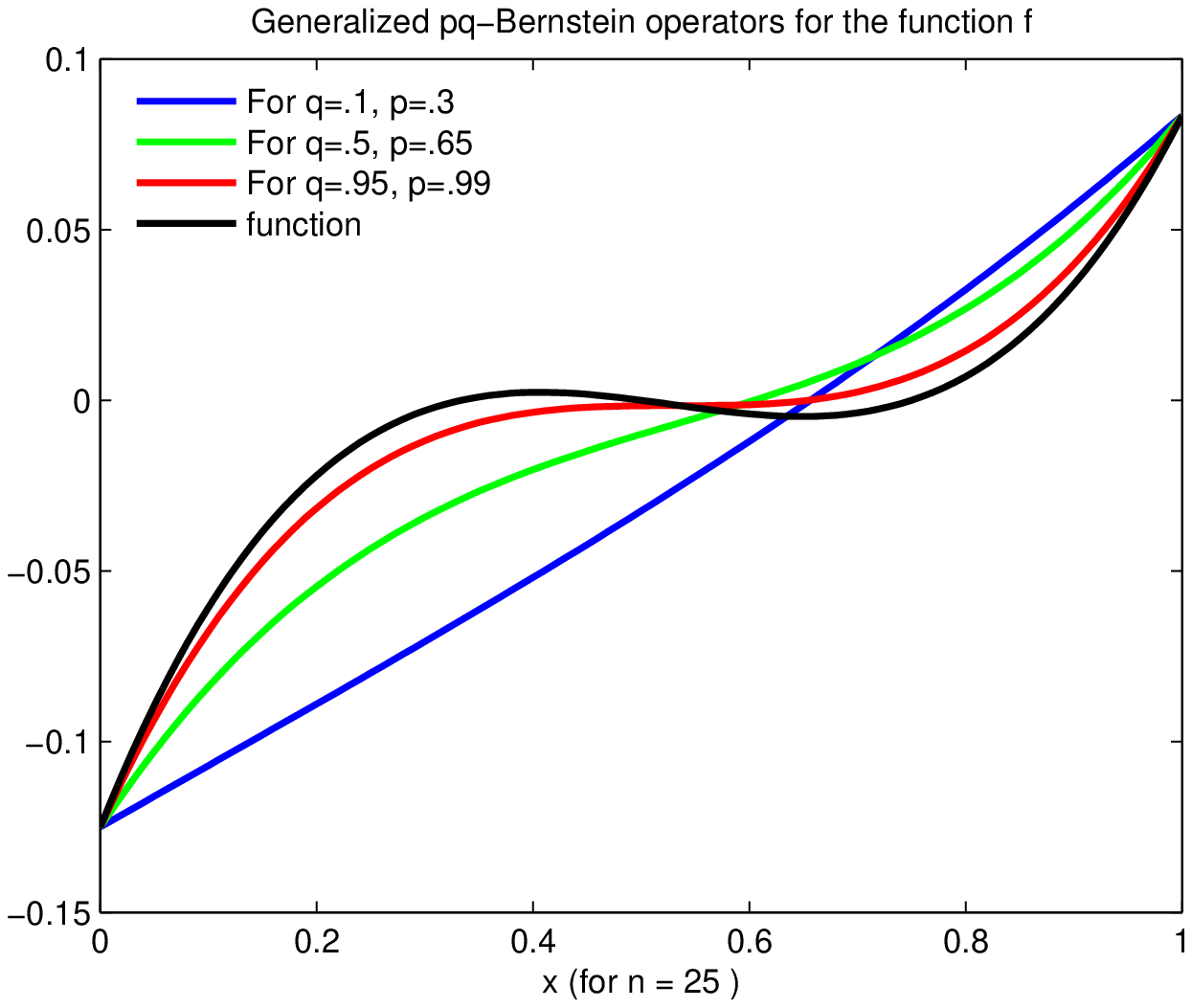}
\end{center}
\caption{}
\end{figure*}

\newpage
In comparison to figure 2, as the value the $n$ increases, operators given by $(2)$ converge towards the function which is shown in figure 3. Also, from figure 2, it can be observed that as the value of $p,q$ approaches towards 1 provided $0<q<p\le1$, operators converge towards the function.
\begin{figure*}[htb!]
\begin{center}
\includegraphics[height=6cm, width=10cm]{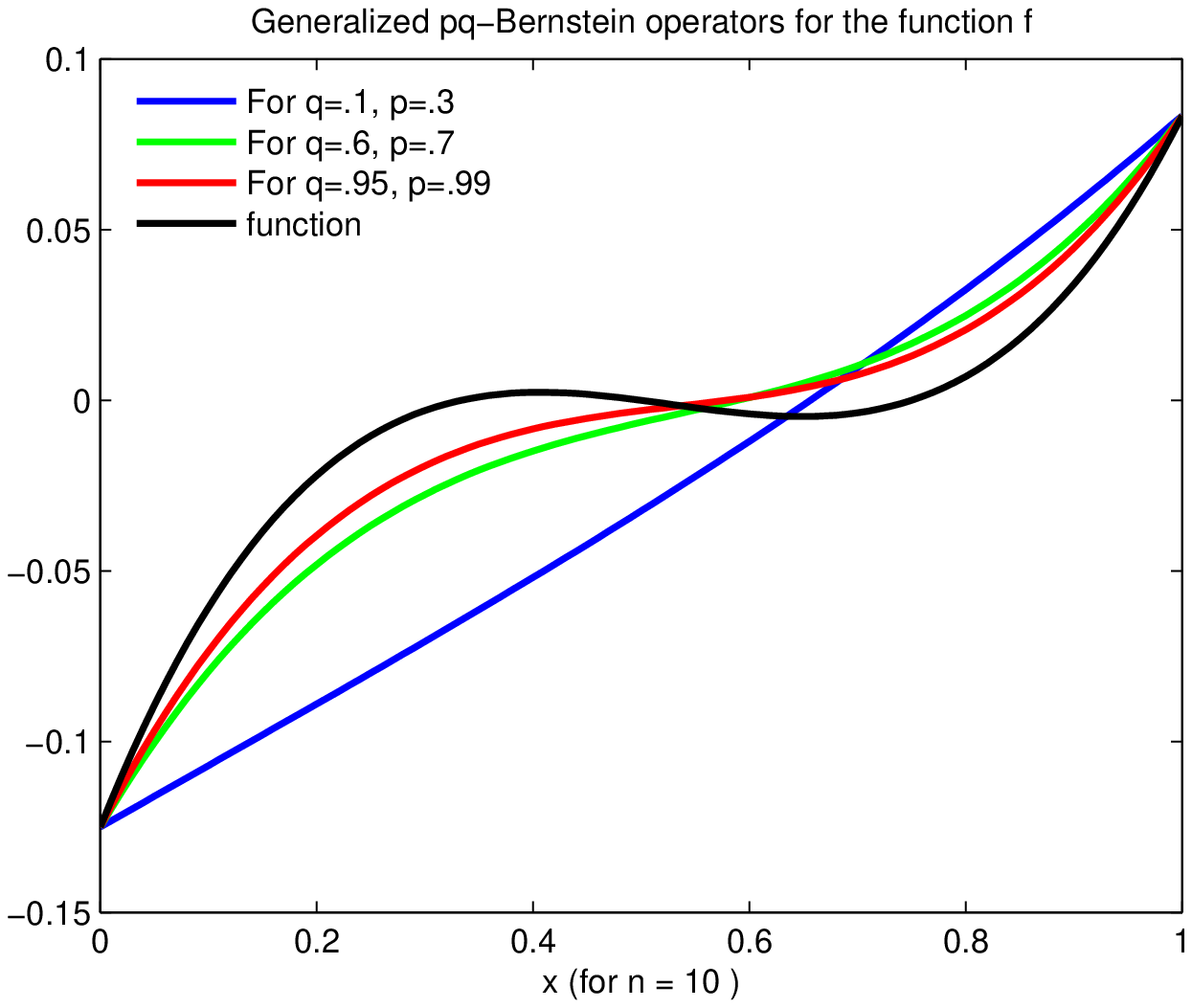}
\end{center}
\caption{}
\end{figure*}

Similarly for different values of parameters $p,q~ \text{and}~ n$
convergence of operators to the function is shown in figure 4 and figure 5.

\begin{figure*}[htb!]
\begin{center}
\includegraphics[height=6cm, width=10cm]{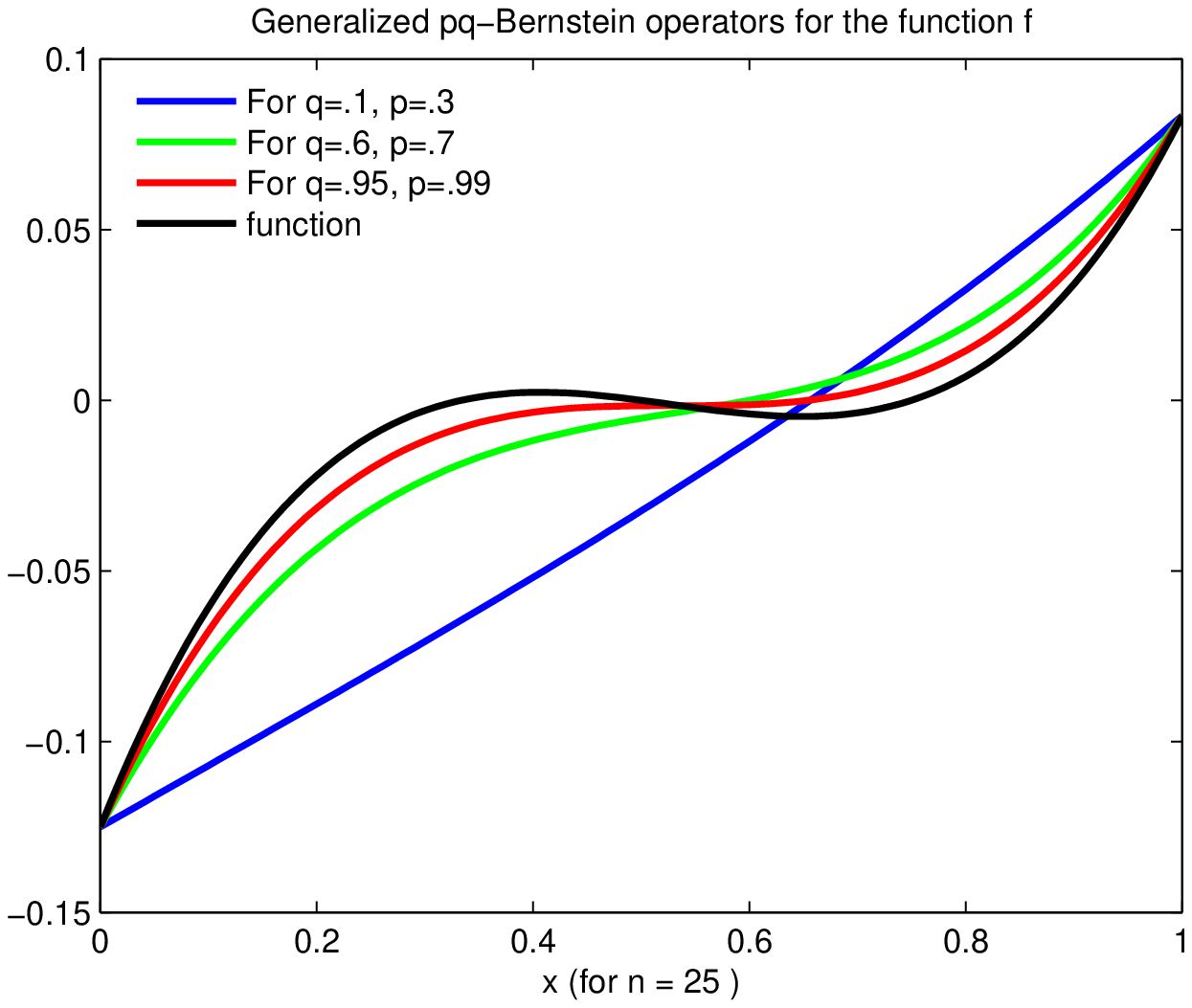}
\end{center}
\caption{}
\end{figure*}

\end{document}